\documentclass{amsart}

\theoremstyle{definition}

\usepackage{hyperref}
\usepackage{breakurl}

\theoremstyle{remark}

\usepackage{algorithm, algpseudocode}

\newcommand{\Tr}{\operatorname{Tr}}

\begin{document}

\title[Algorithm for the characteristic polynomial]{On a fast and nearly division-free algorithm for the characteristic polynomial}


\author{Fredrik Johansson}
\address{Inria Bordeaux-Sud-Ouest and Institut Math.\ Bordeaux,  U. Bordeaux, \newline33400 Talence, France}
\curraddr{}
\email{fredrik.johansson@gmail.com}
\thanks{}



\date{}

\dedicatory{}

\begin{abstract}
We review the Preparata-Sarwate algorithm,
a simple $O(n^{3.5})$ method for computing the
characteristic polynomial, determinant and adjugate of an $n \times n$ matrix
using only ring operations
together with exact divisions by small integers.
The algorithm is a baby-step giant-step version
of the more well-known Faddeev-Leverrier algorithm. We make a few comments
about the algorithm and evaluate its performance empirically.
\end{abstract}

\maketitle

\section{Introduction}

Let $R$ be a commutative ring.
Denote by $\omega > 2$ an exponent of matrix multiplication,
meaning that we can multiply two $n \times n$ matrices
using $O(n^{\omega})$ ring operations (additions, subtractions and multiplications).
Given a matrix $A \in R^{n \times n}$, how fast
can we compute its
characteristic polynomial, determinant and adjugate (or where applicable, inverse), without dividing by elements in $R$?

The obvious division-free algorithm is
cofactor expansion, which uses $O(n!)$ operations.
It is mainly interesting for $n \le 4$ and for sparse symbolic matrices.
The first published efficient method is the Faddeev-Leverrier algorithm~\cite{leverrier1840variations,Helmberg1993,Csanky1976,Hou1998}; Alg.\ 2.2.7 in~\cite{Coh1996}),
which solves the problem using roughly $n$ matrix multiplications,
for a complexity of $O(n^{\omega+1})$.
The Faddeev-Leverrier algorithm requires some divisions,
but it is \emph{nearly division-free} in the sense that all divisors are
fixed small integers $1,2,\ldots,n$ rather than general elements of $R$,
and the divisions are \emph{exact} in the sense that the quotients remain in $R$.

The Berkowitz algorithm \cite{Berkowitz1984} achieves the same complexity
using $O(n^4)$ operations, or $O(n^{\omega+1} \log n)$ with fast
matrix multiplication, without performing any divisions. In practice, the Berkowitz algorithm is faster than the
Faddeev-Leverrier algorithm by a constant factor and
is now widely used in computer algebra systems for linear
algebra over rings where divisions are problematic.

The complexity of the Faddeev-Leverrier algorithm can be improved
using a baby-step giant-step technique,
leading to a method with $O(n^{\omega+0.5}+n^3)$ complexity
which is asymptotically better than the Berkowitz
algorithm.
This method was apparently first discovered by S.~Winograd who did not
publish the result, and then independently discovered and published by Preparata and Sarwate~\cite{Preparata1978}.

Unfortunately, many references to Faddeev-Leverrier algorithm
in the literature present it as an $O(n^4)$ or $O(n^{\omega+1})$ algorithm
without mentioning the improvement of Preparata and Sarwate.
Berkowitz~\cite{Berkowitz1984}
claims that a baby-step giant-step
Faddeev-Leverrier algorithm is possible but only cites private communication with S.~Winograd without giving an explicit algorithm.
The present author is not aware of any software using the Preparata-Sarwate algorithm.

There are perhaps three reasons for the relative obscurity of the method.\footnote{Indeed, the
present author also reinvented the same method. The first version of this
note presented it simply as a ``baby-step giant-step Faddeev-Leverrier algorithm'',
mentioning Berkowitz's attribution of the idea to S.~Winograd,
but without citing the concrete algorithm already given in~\cite{Preparata1978}. The author is grateful to Eric Schost for
pointing out the relevant prior art.}
First of all, the Faddeev-Leverrier algorithm
is numerically unstable~\cite{rehman2011budde,wilkinson1965algebraic}
making it virtually useless for ordinary numerical computation,
and it is in any case much slower than standard $O(n^3)$
methods such as Gaussian elimination and Hessenberg reduction
which are numerically stable.
Second, faster methods are known for the most commonly used exact fields and rings
such as $\mathbb{F}_q$, $\mathbb{Q}$ and $\mathbb{Z}$.
Third, Kaltofen and Villard~\cite{Kaltofen1992,Kaltofen2005,Villard2011} have achieved an even lower
division-free complexity of $O(n^{3.2})$
for the determinant or adjugate and $O(n^{3.29})$ for the characteristic
polynomial with classical multiplication ($\omega = 3$),
and about $O(n^{2.72})$ and $O(n^{2.84})$ for the respective tasks
using the fastest currently known matrix multiplication
algorithm due to Le Gall~\cite{LeGall2014}.
The Kaltofen-Villard algorithm is far more complicated
than the Preparata-Sarwate algorithm, however.

The contribution of this note is twofold. First, we
give explicit pseudocode for a version of the
Preparata-Sarwate algorithm (Algorithm~\ref{alg:flbsgs})
that may serve as a convenient reference for future implementers.
The code makes a superficial change to the algorithm as
it is usually presented~\cite{Preparata1978,Codenotti1990,abdeljaoued2003methodes}, halving
memory usage (this is not an entirely
negligible change
since the Preparata-Sarwate algorithm uses $O(n^{2.5})$ memory).
Second, we discuss the applicability of the algorithm
and perform computational experiments
comparing several algorithms over various rings.

\section{The Faddeev-Leverrier algorithm and its refinement}

We first recall the original Faddeev-Leverrier algorithm (Algorithm~\ref{alg:fl})
for the characteristic polynomial $p_A(x) = c_n x^n + \dots + c_1 x + c_0$ of a matrix $A \in R^{n \times n}$.

\begin{algorithm}[h!]
\caption{Faddeev-Leverrier algorithm}
\label{alg:fl}
\begin{algorithmic}[1]
\Require $A \in R^{n \times n}$ where $n \ge 1$ and $R$ is a commutative ring, $R$ having a unit element and characteristic 0 or characteristic coprime to $1,2,\ldots,n$
\Ensure $(p_A(x), \; \det(A), \; \operatorname{adj}(A))$
\State $c_n = 1, \quad B \gets I, \quad k \gets 1$
\While{$k \le n - 1$}
  \State $B \gets AB$
  \State $c_{n-k} \gets -\frac{1}{k} \Tr(B)$
  \State $B \gets B + c_{n-k} \,I$
  \State $k \gets k + 1$
\EndWhile
\State $c_0 \gets -\frac{1}{n} \Tr(AB)$
\State \Return $(c_0 + c_1 x + \ldots + c_n x^n, \; (-1)^n c_0, \; (-1)^{n+1} B)$
\end{algorithmic}
\end{algorithm}

Algorithm~\ref{alg:fl} is based on the recursion $c_{n-k} = -\frac{1}{k} \sum_{j=1}^k c_{n-k+j} \Tr(A^j)$:
we compute a sequence of matrices (stored in the variable $B$) through repeated
multiplication by $A$, and in each step extract a trace.
The determinant and the adjugate
matrix appear as byproducts of this process: $(-1)^n \det(A)$
is the last coefficient~$c_0$,
and $(-1)^{n+1} \operatorname{adj}(A)$ is the penultimate
entry in the matrix sequence.\footnote{The code can be tightened
assuming $n \ge 2$, in which case the line
before the start of the loop can be changed to
$\{c_n = 1, \; c_{n-1} = -\Tr(A), \; B \gets A + c_{n-1} I, \; k \gets 2\}$.
We can change the loop condition to $k \le n$ and remove the
line after the loop if we omit returning $\operatorname{adj}(A)$.}

It is easy to see that Algorithm~\ref{alg:fl} performs $O(n)$
matrix multiplications and $O(n^2)$ additional arithmetic operations.
The condition on the characteristic of $R$ ensures that we can
divide exactly by each $k$, i.e.\ $(xk)/k = x$ holds for $x \in R, k \le n$.

\begin{algorithm}
\caption{Preparata-Sarwate algorithm (slightly modified)}
\label{alg:flbsgs}
\begin{algorithmic}[1]
\Require $A \in R^{n \times n}$ with the same conditions as in Algorithm~\ref{alg:fl}
\Ensure $(p_A(x), \; \det(A), \; \operatorname{adj}(A))$
\State $m \gets \lfloor \sqrt{n} \rfloor$
\State Precompute the matrices $A^1, A^2, A^3, \ldots, A^m$
\State Precompute the traces $t_k = \Tr(A^k)$ for $k = 1, \ldots, m$
\State $c_n = 1, \quad B \gets I, \quad k \gets 1$
\While{$k \le n-1$}
  \State $m \gets \min(m, n-k)$
  \State $c_{n-k} \gets -\frac{1}{k} \Tr(A, B)$
  \For{$j \gets 1,2,\ldots,m-1$}
    \State $c_{n-k-j} \gets \Tr(A^{j+1}, B)$ \Comment{Using precomputed power of $A$}
    \For{$i \gets 0,1,\ldots,j-1$}
      \State $c_{n-k-j} \gets c_{n-k-j} + t_{j-i} c_{n-k-i}$
    \EndFor
    \State $c_{n-k-j} \gets c_{n-k-j} / (-k-j)$
  \EndFor
  \State $B \gets A^m B$ \Comment{Using precomputed power of $A$}
  \For{$j \gets 0,1,\ldots,m-1$}
    \State $B \gets B + c_{n-k-j} A^{m-j-1}$ \Comment{Using precomputed power, or $A^0 = I$}
  \EndFor
  \State $k \gets k + m$
\EndWhile
\State $c_0 \gets -\frac{1}{n} \Tr(A,B)$
\State \Return $(c_0 + c_1 x + \ldots + c_n x^n, \; (-1)^n c_0, \; (-1)^{n+1} B)$
\end{algorithmic}
\end{algorithm}

Algorithm~\ref{alg:fl} computes a sequence of $n$ matrices
but only extracts a small amount of unique information (the trace) from each
matrix.
In such a scenario, we
can often save time using a baby step giant-step approach in
which we only compute $O(\sqrt{n})$ products explicitly (see \cite{Paterson1973,Brent1978,Berkowitz1984,Johansson2014} for other examples
of this technique).
Preparata and Sarwate~\cite{Preparata1978} improve Algorithm~\ref{alg:fl} by
choosing $m \approx \sqrt{n}$ and precomputing
the powers $A^1,A^2,\ldots,A^m$ and $A^m,A^{2m},A^{3m},\ldots$.
The key observation is that we can compute
$\Tr(A B)$ using $O(n^2)$ operations without forming
the complete matrix product $A B$, by simply evaluating the
dot products for the main diagonal of $AB$.
We denote this \emph{product trace} operation by $\Tr(A, B)$.
Algorithm~\ref{alg:flbsgs} presents a version of this method.
It is clear from inspection that this version performs roughly $m + n/m \approx 2\sqrt{n}$
matrix multiplications of size $n \times n$, and
$O(n^3)$ arithmetic operations in the remaining steps.

Algorithm~\ref{alg:flbsgs} is a small modification of the
Preparata-Sarwate algorithm as it is usually presented.
Instead of computing the giant-step powers $A^m,A^{2m},A^{3m},\ldots$
explicitly, we expand the loop in Algorithm~\ref{alg:fl} to group~$m$ iterations,
using multiplications by $A^m$ to update the running sum over
both the powers and their traces.
This version performs essentially the same number
of operations while using half the amount of memory.

As an observation for implementations, the matrix-matrix multiplications and product traces are done
with invariant operands that get recycled $O(\sqrt{n})$ times.
This can be exploited for preconditioning purposes,
for instance by packing the data more efficiently for arithmetic operations.
We also note that the optimal $m$ may depend on the application,
and a smaller value will reduce memory consumption
at the expense of requiring more matrix multiplications.

\section{Applicability and performance evaluation}

When, if ever, does it make sense to use Algorithm~\ref{alg:flbsgs}?
If $R$ is a field, then the determinant, adjugate and characteristic polynomial
can be computed in $O(n^{\omega})$ operations or $O(n^3)$ classically
allowing divisions~\cite{Dumas2005,Pernet2007,Abdeljaoued2001,neiger2020deterministic}.
Most commonly encountered rings are integral domains,
in which case we can perform computations in the fraction field
and in some cases simply clear denominators.
This does not automatically render an $O(n^{3.5})$ algorithm obsolete,
since naively counting operations may not give the whole picture.
Nevertheless, we can immediately discard some applications:
\begin{itemize}
\item For computing over $\mathbb{R}$ and $\mathbb{C}$ in ordinary floating-point arithmetic,
the Faddeev-Leverrier algorithm is slower and far less numerically stable than textbook
techniques such as reduction to Hessenberg form and Gaussian elimination with $O(n^3)$ or better complexity, as we already noted in the introduction.
\item For finite fields, classical $O(n^3)$ methods using divisions have no drawbacks, and linear algebra with $O(n^{2.81})$ Strassen complexity is well
established~\cite{dumas2012computational}. Over rings with small characteristic, the applicability of Algorithm~\ref{alg:flbsgs} is in any case limited due to the integer divisions.
\end{itemize}

Generally speaking, division-free or nearly division-free algorithms are interesting
for rings and fields $R$ where dividing recklessly can lead to
coefficient explosion (for example, $\mathbb{Q}$)
or in which testing for zero is problematic (for example,
exact models of~$\mathbb{R}$).
The optimal approach in such situations is usually to avoid computing directly in $R$ or its fraction field,
for example using modular arithmetic and interpolation techniques,
but such indirect methods are more difficult to implement
and must typically be designed on a case by case basis.
By contrast, Algorithm~\ref{alg:flbsgs} is easy to use anywhere.
We will now look at some implementation experiments.

\subsection{Integers}

For exact linear algebra over $\mathbb{Z}$ and~$\mathbb{Q}$,
the best methods are generally fraction-free versions of
classical algorithms (such as the Bareiss version of Gaussian
elimination) for small $n$, and multimodular or $p$-adic
methods for large~$n$ (see for example~\cite{Dumas2005,Pernet2007}).
We do not expect Algorithm~\ref{alg:flbsgs} to beat
those algorithms, but it is instructive to examine its performance.
Table~\ref{table:timingsint} shows timings for
computing a determinant, inverse or characteristic polynomial of
an $n \times n$ matrix over~$\mathbb{Z}$ with random entries in $-10,\ldots,10$,
using the following algorithms:

\begin{itemize}
\item FFLU: fraction-free LU factorization using the Bareiss algorithm.
\item FFLU2: as above, but using the resulting decomposition to compute $A^{-1}$ (equivalently determining $\operatorname{adj}(A)$) by solving $AA^{-1} = I$.
\item ModDet: a multimodular algorithm for the determinant.
\item ModInv: a multimodular algorithm for the inverse matrix.
\item ModCP: a multimodular algorithm for the characteristic polynomial.
\item Berk: the Berkowitz algorithm for the characteristic polynomial.
\item Alg1: original Faddeev-Leverrier algorithm, Algorithm~\ref{alg:fl}.
\item Alg2: modified Preparata-Sarwate algorithm, Algorithm~\ref{alg:flbsgs}.
\end{itemize}

We implemented Alg1 and Alg2 on top of Flint~\cite{Har2010},
while all the other tested algorithms are builtin Flint methods.

\begin{table}
\begin{center}
\caption{Time in seconds to compute characteristic polynomial (C), determinant (D), adjugate/inverse (A) of an $n \times n$ matrix over $\mathbb{Z}$ with random elements in $-10,\ldots,10$, using various algorithms.}
\label{table:timingsint}
\renewcommand{\baselinestretch}{1.15}
\begin{footnotesize}
\begin{tabular}{ c | c c | c c | c | c c c } 
$n$  & FFLU       & ModDet & FFLU2  & ModInv   & ModCP & Berk & Alg1 & Alg2 \\
     &    D       &   D      &  DA      &   A     & CD  &  CD  & CDA   & CDA \\ \hline
10   & 0.0000060  & 0.000021 & 0.000015 & 0.00012  & 0.000016 & 0.000035 & 0.000015 & 0.000030 \\
20   & 0.000036   & 0.000078 & 0.00043  & 0.00096  & 0.00011  & 0.0010   & 0.00086  & 0.00061 \\
50   & 0.0023     & 0.0012   & 0.011    & 0.016     & 0.0039   & 0.048    & 0.052     & 0.017 \\
100  & 0.039      & 0.0068   & 0.14     & 0.18     & 0.055    & 0.84     & 1.1       & 0.22 \\
200  & 0.64       & 0.044    & 2.3      & 1.8     & 0.89     & 16       & 27        & 4.1 \\
300  & 3.4        & 0.15     & 13       & 9.0      & 4.6      & 94       & 174       & 20 \\
400  & 12         & 0.38     & 45       & 22      & 15       & 321      & 696       & 66 \\
500  & 32         & 0.77     & 127      & 52     & 37       & 900      & 2057      & 150 \\
\end{tabular}
\end{footnotesize}
\end{center}
\end{table}

\subsubsection{Observations}

Alg2 clearly outperforms both Alg1 and Berk for large $n$,
making it the best algorithm
for computing the characteristic
using direct arithmetic in $\mathbb{Z}$
(the modular algorithm is, as expected, superior).
Alg2 is reasonably competitive for computing the inverse
or adjugate matrix, coming within a factor 2-3$\times$ of FFLU2
and ModInv in this example.
For determinants, the gap to the FFLU algorithm is
larger, and the modular determinant algorithm
is unmatched.

\subsection{Number fields}

Exact linear algebra over algebraic number fields $\mathbb{Q}(a)$
is an interesting use case for division-free algorithms
since coefficient explosion is a significant problem
for classical $O(n^3)$ algorithms.
As in the case of $\mathbb{Z}$ and $\mathbb{Q}$, modular
algorithms are asymptotically more efficient
than working over $\mathbb{Q}(a)$ directly, but harder
to implement. Here we compare the following algorithms:

\begin{itemize}
\item Sage: the \texttt{charpoly} method in SageMath~\cite{Sag2020}, which implements a special-purpose algorithm for cyclotomic fields based on modular computations and reconstruction using the Chinese remainder theorem.
\item Hess: Hessenberg reduction for the characteristic polynomial
\item Dani: Danilevsky's algorithm for the characteristic polynomial.
\item LU: LU factorization to compute the determinant.
\item FFLU: fraction-free LU factorization using the Bareiss algorithm.
\item LU2 and FFLU2: as above, but using the resulting decomposition to compute $A^{-1}$ (equivalently determining $\operatorname{adj}(A)$) by solving $AA^{-1} = I$.
\item Berk (Berkowitz), Alg1 and Alg2 as in the previous section.
\end{itemize}

With the exception of the Sage function,
we implemented all the algorithms
using Antic~\cite{Har2015} for number field arithmetic
and Flint for other operations.
We perform fast matrix multiplication
by packing number field elements into integers and multiplying matrices over $\mathbb{Z}$ via Flint.
Our implementations of LU, LU2, Alg1 and Alg2 benefit from matrix
multiplication while Hess, Dani, FFLU, FFLU2 and Berk do not.
The benchmark is therefore not representative
of the performance that ideally should be achievable with these algorithms,
although it is fair in the sense
that the implementation effort for Alg1 and Alg2 was minimal
while the other algorithms would require much more code to speed up
using block strategies.

Table~\ref{table:timingscyclo} compares timings for two kinds of input:
random matrices over a fixed cyclotomic field,
and discrete Fourier transform (DFT) matrices which have special structure.
Choosing cyclotomic fields allows us to compare with the dedicated
algorithm for characteristic polynomials in Sage;
the corresponding method for generic number fields in Sage is far slower.
All the other algorithms make no assumptions about the field.

\begin{table}
\begin{center}
\caption{Time in seconds to compute characteristic polynomial (C), determinant (D), adjugate/inverse (A) of a matrix over a cyclotomic number field.}
\label{table:timingscyclo}
\renewcommand{\baselinestretch}{1.15}
\begin{footnotesize}
\begin{tabular}{ c | c | c c | c c | c c | c c c } 
$n$ & Sage & \!Hess\! & \!\!\!Dani\!\! & \!LU\! & \!FFLU\! & \!LU2\! & \!FFLU2\! & \!Berk\!\! &  \!\!\!Alg1\!\!\! & \!\!\!Alg2\!\!\! \\
    & CD & CD &  CD   & \!D        & D  & DA    & \!DA  & CD &   \!CDA\!\!   & \!\!CDA\!\!   \\
\hline
\multicolumn{11}{c}{Input: $n \times n$ matrix over $\mathbb{Q}(\zeta_{20})$, entries $\sum_k (p/q) \zeta_{20}^k$, random $|p| \le 10$, $1 \le q \le 10$. } \\
\hline
10 & 0.038  & 0.31 & 0.16 & 0.024 & 0.0059 & 0.21 & 0.11        & 0.010 & 0.0073 & 0.010 \\
20 & 0.12   & 19 & 6.7 & 0.22 & 0.067 & 2.6 & 1.4               & 0.28  &  0.15 & 0.16 \\
30 & 0.39   & 200 & 67 & 0.93 & 0.31 & 12 & 6.8                 & 2.0   &  1.1 & 0.8 \\
40 & 1.1    &    & 353 & 2.8 & 0.9 & 37 & 22                     & 7.5   & 3.7 & 2.6 \\
50 & 1.9    &    &      & 7.0 & 2.3 & 88 & 56                   & 22    & 8.7 & 5.7 \\
60 & 3.4    &    &      & 15 & 4.7 & 182 & 119                  & 54    & 19  & 12 \\
70 & 5.1    &    &     & 29 & 8.6 & 337 & 230                   & 114   & 39  & 22 \\
80 & 7.5    &    &     & 53 & 15 & 581 & 409                    & 208   & 67 & 34 \\
90 & 11     &    &     & 89 & 24 & &                            & 397   & 144 & 54 \\
100 & 15    &    &      & 144 & 41 & &                           & 608  & 670 & 130 \\
120 & 24    &    &      & 322 & 83 & &                           & 1439 & 3013 & 420 \\
\hline
\multicolumn{11}{c}{Input: $n \times n$ DFT matrix over $\mathbb{Q}(\zeta_{n})$, entries $A_{i,j} = \zeta_n^{(i-1)(j-1)}$.  } \\
\hline
10 & 0.010  & 0.0018 & 0.0016 & \!\!0.00017\!\! & \!\!0.00022\!\! & \!\!0.00076\!\! & 0.0014 & \!\!0.00061\!\! & \!\!0.00075\!\! & \!\!0.00059\!\! \\
20 & 0.039  & 0.0019 & 0.0024 & 0.0017  & 0.0046  & 0.0071  & 0.038 & 0.020 & 0.020 & 0.0070 \\
50 & 1.3  & 0.17   & 0.13   & 0.065   & 0.80    & 0.28    & 6.0   & 8.2   & 2.0   & 0.49   \\
100 & 22   & 5.4    & 22     & 0.89    & 43      & 5.3     & 335   & 803   & 223   & 29 \\
150 & 78   & 22     & 7.9    & 4.4     & 214     & 19      & 1423  & 7259   & 933   & 138 \\
200 & 333  & 1928       & 140       &  31       &    1655     & 192        &       &       &       &  1687    \\ \hline
\end{tabular}
\end{footnotesize}
\end{center}
\end{table}

\subsubsection{Observations}

There are no clear winners since there is a
complex interplay between
operation count,
multiplication algorithms,
matrix structure and coefficient growth.
Modular algorithms are the best solution
in general for large $n$, but implementations for number
fields are complex and less readily available in
current software than for $\mathbb{Z}$ and $\mathbb{Q}$.

Among the non-modular algorithms,
the $O(n^3)$ division-heavy Hessenberg and Danilevsky algorithms
are nearly useless due to coefficient explosion for generic
input, but both perform well on the DFT matrix.
The LU and FFLU algorithms have more even performance
but alternate with each other for the advantage depending on the matrix.
Alg2 has excellent average performance for the determinant,
characteristic polynomial as well as the adjugate matrix
considering the large
variability between the algorithms for different input.
It is highly competitive for computing the inverse or adjugate
of the random matrix.

\subsection{Ball arithmetic}

Division-free algorithms are useful when computing
rigorously over $\mathbb{R}$ and $\mathbb{C}$ in interval arithmetic
or ball arithmetic. The reason is that we cannot test whether elements
are zero, so algorithms like Gaussian elimination and Hessenberg
reduction fail when they need to branch upon
zero pivot elements
or zero vectors.
Although zeros will not occur for random input,
they are likely to occur for structured matrices arising in applications.
\emph{A posteriori} verification of approximate numerical solutions
or perturbation analysis is in principle the best workaround~\cite{Rump2010},
but it is sometimes useful to fall back to more direct division-free methods,
especially when working in very high precision.

\begin{table}
\begin{center}
\caption{Time in seconds to compute characteristic polynomial (C), determinant (D), adjugate/inverse (A) of a random $n \times n$ matrix in real ball arithmetic. The respective algorithms were run with $333 + p$ bits of precision, with $p$ chosen to give roughly 100-digit output accuracy.}
\label{table:timingsarb}
\renewcommand{\baselinestretch}{1.15}
\begin{footnotesize}
\begin{tabular}{ c | c c c | c | c c | c c c } 
$n$  & Hess & Hess2 & Dani & Eig & LU & LU2 & Berk & Alg1 & Alg2 \\
     & CD    & CD    &  CD  & CD & D  & D      & CD  & CDA   & CDA   \\ \hline
10 & 0.00021 & 0.00038 &  0.00022    & 0.017 & 0.000068 & 0.00023 & 0.00080 & 0.0010 & 0.00078 \\
20 & 0.0021  & 0.0032 &  0.0020    & 0.18  & 0.00052  & 0.0015  & 0.0039 & 0.017  & 0.0092  \\
50 & 0.045   & 0.057 &  0.048    & 4.6   & 0.0078   & 0.019   & 0.22   & 0.61   & 0.21    \\
100 & 0.64   & 0.69 &  0.61    & 56    & 0.062    & 0.15    & 6.3    & 9.7    & 2.4     \\
150 & 3.5    & 3.5 &  3.0    & 245  & 0.23     & 0.44    & 52     & 52     & 11      \\
200 & 12     & 11  &  10    &       & 0.59     & 1.0     & 224    & 176    & 34      \\
250 & 31     & 29  &  25    &       & 1.4      & 1.9     & 687    & 460    & 73      \\
300 & 66     & 59  &  53    &       & 2.5      & 3.2     & 1804   & 1075    & 160     \\
350 & 135    & 115 &  110    &       & 4.4      & 5.0     & 4033   & 2107   & 306     \\ \hline
$p$      & $10n$ & $6n$  &  $10n$ &  $0$  &  $n$  & $0$  &  $6n$ & $6n$ & $6n$ \\
\end{tabular}
\end{footnotesize}
\end{center}
\end{table}

Table~\ref{table:timingsarb} shows timings for
computing a determinant or characteristic polynomial with 100-digit
accuracy using the following algorithms implemented in ball arithmetic.
The input is taken to be an $n \times n$ real matrix with uniformly random entries in $[0,1]$.
For this experiment,
we only focus on the determinant and characteristic polynomial
(the conclusions regarding matrix inversion would be similar to those regarding the determinant).

\begin{itemize}
\item Hess: Hessenberg reduction using Gaussian elimination.
\item Hess2: Hessenberg reduction using Householder reflections.
\item LU: LU factorization using Gaussian elimination.
\item LU2: approximate computation of the determinant using LU factorization followed by \emph{a posteriori} verification.
\item Eig: approximate computation of the eigenvalues using the QR algorithm followed by \emph{a posteriori} verification and reconstruction of the characteristic polynomial from its roots.
\item Berk (Berkowitz), Alg1 and Alg2 as in the previous section.
\end{itemize}

All algorithms were implemented in Arb\cite{Joh2017} which uses
the accelerated dot product and matrix
multiplication algorithms described in~\cite{Johansson2019}.
The LU, LU2, Alg1 and Alg2 implementations benefit
from fast matrix multiplication while Hess, Hess2, Eig and Berk do not.

The methods LU2 and Eig are numerically stable: the output
balls are precise to nearly full precision for well-conditioned input.
All other algorithms are unstable in ball arithmetic
and lose $O(n)$ digits of accuracy.
At least on this example, the rate of loss is almost the same
for Hess2, Berk, Alg1 and Alg2, while LU is more stable
and Hess and Dani are less stable.
To make the comparison meaningful, we set the working precision
(shown in the table) to an experimentally determined
value so that all algorithms enclose
the determinant with around 100 digits of accuracy.

\subsubsection{Observations}

For computing the characteristic polynomial
in high-precision ball arithmetic, it seems prudent to try
Hessenberg reduction and fall back to a division-free algorithm when it fails
due to encountering a zero vector.
The Berkowitz algorithm is the best fallback for small $n$,
while Alg2 wins for large $n$ ($n \approx 50$, although the
cutoff will vary).
On this particular benchmark, Alg2 runs only
about $4\times$ slower than
Hessenberg reduction,
making it an interesting \emph{one-size-fits-all} algorithm.
The verification method (Eig) gives the best results if the
precision is constrained, but is far more expensive than the other methods.

For computing the determinant alone,
all the division-free methods are clearly
inferior to methods based on LU factorization
in this setting. The only
advantage of the division-free algorithms is that
they are foolproof while LU factorization
requires some attention to implement correctly.

Better methods for computing the characteristic polynomial
in ball arithmetic or interval arithmetic are surely possible.
For the analogous problem of computing the characteristic polynomial
over $\mathbb{Q}_p$, see~\cite{Caruso2017}.

\subsection{Polynomial quotient rings}

At first glance Algorithm~\ref{alg:flbsgs} seems to hold
potential for working over multivariate polynomial quotient rings.
Such rings need not be integral domains
and division can be very expensive (requiring Gr\"{o}bner basis computations).
Unfortunately, in most examples we have tried, Algorithm~\ref{alg:flbsgs}
performs worse than both the Berkowitz algorithm and
Algorithm~\ref{alg:fl},
presumably because repeated multiplication
by the initial matrix $A$ is much cheaper
than multiplication by a power of $A$ which generally
will have much larger entries.
There may be special classes of matrices for which the
method performs well, however.

\section{Discussion}

We find that the Preparata-Sarwate algorithm often
outperforms the Berkowitz algorithm in practice,
in some circumstances even being competitive with $O(n^3)$ algorithms.
We can therefore recommend it for more widespread adoption.

Galil and Pan \cite{Galil1989,abdeljaoued2003methodes}
have further refined the
Preparata-Sarwate algorithm to eliminate the $O(n^3)$
complexity term which dominates asymptotically with a sufficiently fast
matrix multiplication algorithm.
We have not tested this method since
those $O(n^3)$ operations are negligible in practice.

An interesting problem is whether it is possible to design a division-free
or nearly division-free algorithm with
better than $O(n^4)$ classical complexity that minimizes
the observed problems with growing entries, particularly in multivariate rings.

\bibliographystyle{plain}
\bibliography{references}

\end{document}